\documentclass[a4paper,12pt]{amsart}
\usepackage{amscd,amsfonts,amssymb}
\usepackage{graphicx,tikz}
\usepackage{mathrsfs}
\usetikzlibrary{arrows}
\usepackage{enumerate}
\usepackage[shortlabels]{enumitem}
\usepackage{mathrsfs}
\usepackage{amssymb,amsmath,amsthm,color}
\usepackage{caption,subcaption}
\usepackage{hyperref}
\usepackage{url}
\usepackage{setspace}
\usepackage{float}

\date{\today}

\allowdisplaybreaks

\newtheorem{theorem}{Theorem}[section]

\newtheorem{lemma}[theorem]{Lemma}

\newtheorem{remark}[theorem]{Remark}

\def\be#1 {\begin{equation} \label{#1}}
\newcommand{\ee}{\end{equation}}

\def\sqw{\hbox{\rlap{\leavevmode\raise.3ex\hbox{$\sqcap$}}$%
\sqcup$}}
\def\findem{\ifmmode\sqw\else{\ifhmode\unskip\fi\nobreak\hfil
\penalty50\hskip1em\null\nobreak\hfil\sqw
\parfillskip=0pt\finalhyphendemerits=0\endgraf}\fi}

\newcommand{\R}{{\mathbb {R}}}

\author{Saurabh Shrivastava and Kalachand Shuin}
\address{
}
\email{}
\address[Saurabh Shrivastava and Kalachand Shuin]{
Department of Mathematics\\
Indian Institute Science Education and Research Bhopal\\
Bhopal-462066, India}
\email{\{saurabhk,kalachand16\}@iiserb.ac.in}


\begin{document}

\title[multilinear convolution operators]{$L^p$ estimates for multilinear convolution operators defined with spherical measure}

\begin{abstract}
Let $\sigma=(\sigma_{1},\sigma_{2},\dots,\sigma_{n})\in \mathbb{S}^{n-1}$ and $d\sigma$ denote the normalised Lebesgue measure on $\mathbb{S}^{n-1},~n\geq 2$. For functions $f_1, f_2,\dots,f_n$ defined on $\R$ consider the multilinear operator given by
$$T(f_{1},f_{2},\dots,f_{n})(x)=\int_{\mathbb{S}^{n-1}}\prod^{n}_{j=1}f_{j}(x-\sigma_j)d\sigma, ~x\in \R.$$
In this paper we obtain necessary and sufficient conditions on exponents $p_1,p_2,\dots,p_n$ and $r$ for which the operator $T$ is bounded from $\prod_{j=1}^n L^{p_j}(\R)\rightarrow L^r(\R),$ where $1\leq p_j,r\leq \infty, j=1,2,\dots,n.$ This generalizes the results obtained in~\cite{jbak,oberlin}. 
\end{abstract}
\maketitle

\section{Introduction and preliminaries}
Let $\sigma=(\sigma_{1},\sigma_{2},\dots,\sigma_{n})$ denote a point on the unit sphere $\mathbb{S}^{n-1}\subset \R^n,~n\geq 2,$ and $d\sigma$ be the normalised Lebesgue measure on $\mathbb{S}^{n-1}.$ For continuous functions $f_1, f_2,\dots, f_n$ on $\R$ consider the multilinear spherical averages defined by 
\begin{eqnarray}\label{def}T(f_{1},f_{2},\dots,f_{n})(x)=\int_{\mathbb{S}^{n-1}}\prod^{n}_{j=1}f_{j}(x-\sigma_j)d\sigma,~x\in \R.
\end{eqnarray}
Let $1\leq p_1,p_2,\dots,p_n,r\leq \infty$. We are interested in studying the $L^p$ estimates for the operator $T$ at $(p_1,p_2,\dots,p_n;r)$, i.e., 
\begin{eqnarray}\label{bdd}
\|T(f_{1},f_{2},\dots,f_{n})\|_r\lesssim \prod_{j=1}^n \|f_j\|_{p_j}.
\end{eqnarray}
 
The notation $A\lesssim B$ (and $A\gtrsim B$) means that there exists an implicit constant $C>0,$ such that $A\leq CB$ (and $A\geq CB$). We will not keep track of the constants and often use the notation as mentioned above. We will also require weaker notion of boundedness of operators between Lorentz spaces as the operator may not always satisfy strong type estimates. 

We need to consider a general form of the operator $T$ as it would be required in many of our proofs. Let $\{v_{1},v_{2},\dots ,v_{n}\}$ be linearly independent vectors in $\R^n$. Consider the following general form of the operator $T$ given by  $$T_v(f_{1},\dots,f_{n}):=\int_{\mathbb{S}^{n-1}}\prod^{n}_{j=1}f_{j}(x-v_{j}\cdot\sigma)d\sigma.$$
Recall that for a given multilinear ($n-$linear) operator one can consider $n$ adjoint operators associated with it. More specificly, we have adjoints of $T_v$ given by 
$$\langle T_v^{*j}(f_1,f_2,\dots,f_n),h\rangle :=\langle T_v(f_1,f_2,\dots,f_{j-1},h,f_{j+1},\dots, f_n),f_j\rangle.$$
It is easy to verify that $T_v^{*j}$ is similar to $T_v$  with a different set of linearly independent vectors than that of $T_v$. Using duality arguments boundedness of $T_v$ at $(p_1,p_2,\dots,p_n;r)$ implies the corresponding result for $T_v^{*j}$ at $(p_1,p_2,\dots,p_{j-1},r',p_{j+1},\dots,p_{n};p'_j)$. We will refer to these points as dual points to each other. Here $p'$ denotes the conjugate index to $p$ given by $\frac{1}{p}+\frac{1}{p'}=1$.
 
In~\cite{oberlin} Oberlin established nesessary and sufficent conditions for the boundedness of $T$ from $\prod_{j=1}^nL^{p}(\R)\rightarrow L^r(\R)$. Later, Bak and Shim~\cite{jbak} extended Oberlin's result improving the range of $p$ and $r$ for the strong type boundedness of $T$. We also refer to ~\cite{oberlin1} for Young's inequality for multilinear convolution operators. In order to describe the known results we require some notation.  

Let $R=R(n)$ denote the closed convex hull in $\R^2$ of points $O=(0,0), B=(\frac{n-1}{n+1},0), M=(\frac{n+1}{n+3},\frac{2}{n+3}), A=(\frac{n+1}{n+2},1), F=(\frac{1}{n},1),$ see Figure~$1$ for detail. Note that the point $(\frac{1}{p},\frac{1}{r})$ corresponds to $(p,p,\dots,p;r)$. In~\cite{oberlin} Oberlin proved the following result concerning the boundedness of the operator $T$.
\begin{theorem}\cite{oberlin}\label{oberlin} If the operator $T$ is of strong type at $(p,p,\dots,p;r)$ then $(\frac{1}{p},\frac{1}{r})$ lies in the region $R$.  Conversely, if $(\frac{1}{p},\frac{1}{r})$ lies in the region $R$ and not on the two closed line segments $AM$ and $MB$, then  $T$ is of strong type at $(p,p,\dots,p;r)$. Further, for points $(\frac{1}{p},\frac{1}{r})$ lying on the line segments $AM$ and $MB$, the operator $T$ is of restricted type at $(p,p,\dots,p;r)$, i.e., estimate ~(\ref{bdd}) holds at $(p_1,p_2,\dots,p_n;r)$ for $f_j's$ restricted to characteristic functions. 
\end{theorem}

\begin{figure}[H]
	\includegraphics[width=.5\textwidth]{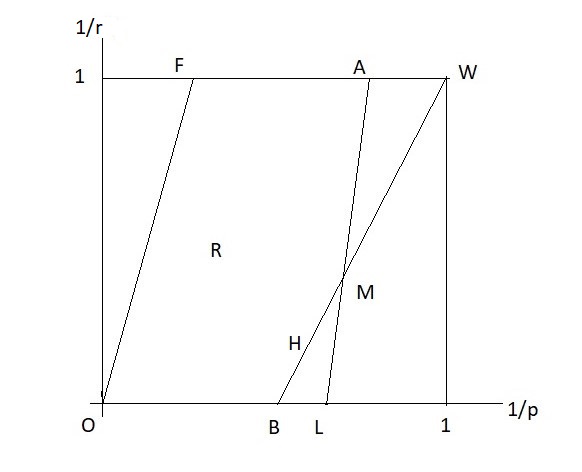}
	\vspace{0mm}
	\caption{Region R}	
\end{figure}
The question of strong type boundedness of the operator $T$ at points lying on the closed line segments $AM$ and $MB$ remained unresolved for a long time. In 1998, Bak and Shim~\cite{jbak} settled this question and filled in the gap between necessary and sufficient conditions in Theorem~\ref{oberlin} for dimension $n\geq 3.$  More precisely, they proved the following. 

\begin{theorem}\cite{jbak}\label{jbak} For $n\geq 3$, the operator $T$ is of strong type at $(p,p,\dots,p;r)$ if, and only if $(\frac{1}{p},\frac{1}{r})$ lies in the region $R$.
\end{theorem} 
Bak and Shim~\cite{jbak} also addressed the question in dimension $n=2$. They obtained the following positive and negative results in this case. 
\begin{theorem}\cite{jbak}\label{jbak1} 
In dimension $n=2$, the following results hold. 
\begin{enumerate}
\item $T$ is of strong type at $(\frac{4}{3},\frac{4}{3};1)$.
\item $T$ is bounded from $L^{3,s}(\R)\times L^{3,s}(\R)\rightarrow L^{\infty}(\R)$
if, and only if $0<s\leq 2$. In particular, $T$ fails to be of strong type at $(3,3;\infty)$. 
\item Let $H$ denote the point $(\frac{1}{2},\frac{1}{4})$ on the line segment $BM$. If $(\frac{1}{p},\frac{1}{r})$ lies on either of the closed line segments $AM$ and  $MH$ then $T$ is of strong type at $(p,p;r)$.
\end{enumerate}
\end{theorem}
Further, we note that in dimension $n=2$ the operator $T$ coincides with the bilinear spherical averaging operator. The bilinear spherical averages and the corresponding bilinear spherical maximal function have been studied by several authors in the recent past. 
For functions $f,g \in \mathcal{S}(\R^d), d\geq 1$, the bilinear spherical average is defined by
\begin{eqnarray*}
\mathcal{A}(f,g)(x):=\int_{\mathbb{S}^{2d-1}}f(x-y)g(x-z)d\sigma(y,z).
\end{eqnarray*}
where $d\sigma(y,z)$ is the normalised Lebesgue measure on the sphere $\mathbb{S}^{2d-1}$. 

Observe that the operator $\mathcal A$ for $d=1$ is same as $T$ for $n=2$. The operator $\mathcal A$ and the corrsponding bilinear maximal function was introduced and studied in~\cite{GGIP}. Later, in ~\cite{BGHH, GHH} authors estbalished partial results obtaining $L^{p_1}(\R^d)\times L^{p_2}(\R^d)\rightarrow L^p(\R^d)$ estimates for the bilinear spherical maximal operator for a certain range of $p_1, p_2$ and $p$ with some assumptions on the dimension $d$. Very recently, in~\cite{JL} Jeong and Lee proved $L^{p_1}(\R^d)\times L^{p_2}(\R^d)\rightarrow L^p(\R^d)$ estimates for the maximal operator for the best possible range of exponents $p_1, p_2$ and $p$ for all $d\geq 2.$ They also obtained $L^p$ improving estimates for the bilinear spherical averaging operator $\mathcal A$ for $d\geq 2$. However, the case of dimension $d=1$ has not been addressed so far. We shall fill this gap in this paper. We also refer to the recent papers~\cite{AP,Do, LSK} for further generalisation of the bilinear spherical maximal functions to the multilinear and product type setting. 

In this paper our aim is to establish necessary and sufficient conditions on exponents $1\leq p_j, r\leq \infty, j=1,2,\dots,n$ for the $(p_1,p_2,\dots,p_n;r)$ boundedness of the operator $T$. 
Note that the results due to Oberlin~\cite{oberlin} and Bak and Shim~\cite{jbak} addressed the same when $p_j=p$ for all $j=1,2,\dots,n.$ We extend their results to the full possible range of exponents and thereby allow the possibility of $p_j$ assuming different values. Moreover, as pointed out earlier the boundedness of $T$ for $n=2$ yields the corresponding $L^p$ improving estimates for the bilinear averaging operator $\mathcal A$ for $d=1$.  Our proofs are motivated from the ideas presented in~\cite{jbak,oberlin}. Along with the standard multilinear interpolation theorems the following multilinear interpolation result due to Christ~\cite{Ch} plays an important role.
\begin{lemma}\label{christ}\cite{Ch}
Let $n\geq 2$ and $\Sigma$ be a nontrivial closed $(n-1)$-simplex in the unit cube $[0,1]^n$. Assume that the hyperplane $\Gamma$ generated by $\Sigma$ is not parallel to any of the coordinate axes. If $S$ is a multilinear ($n-$linear) operator such that it is bounded from $L^{p_1,1} \times \dots\times L^{p_n,1}\rightarrow Y$ at all endpoints $(\frac{1}{p_1},\frac{1}{p_2},\dots,\frac{1}{p_n})$ of the simplex $\Sigma$, where $Y$ is a Banach space. Then for $(\frac{1}{p_1},\frac{1}{p_2},\dots,\frac{1}{p_n})$ an interior point of $\Sigma$ and $1\leq q_i\leq \infty$ satisfying $\sum\limits_{i=1}^n \frac{1}{q_i}=1$, the operator $S$ is bounded from  $L^{p_1,q_1}\times \dots\times L^{p_n,q_n}\rightarrow Y.$
\end{lemma} 

The remaining part of the paper is organized as follow. In Section~\ref{main} we state the main results and describe the necessary region. In Section~\ref{proof:nec} we establish the necessary conditions on exponents for boundedness of the operator. Section~\ref{proof:n=2} is devoted to proving boundedness of the operator for points in the necessary region for $n=2$. Finally, in Section~\ref{proof:all} we complete the proof for $n\geq3.$
\section{Main results}\label{main}
\subsection{Necessary part}
The following result describes necessary conditions on exponents for $L^p$ boundedness of the operator $T$. 
\begin{theorem}\label{mainresult:nec}{\bf Necessary conditions.}  Let $1\leq p_j,r\leq \infty, j=1,2,\dots,n,$ be given exponents. If the operator $T$ is of strong type at $(p_1,p_2,\dots,p_n;r)$ then the following conditions hold.
\begin{enumerate}[i)]
	\item  $\frac{1}{r}\leq\sum^{n}_{j=1}\frac{1}{p_{j}}$, 
	\item $\sum^{n}_{j=1,j\neq k}\frac{1}{p_{j}}+\frac{2}{p_{k}}\leq n-1+\frac{2}{r},$ for $1\leq k\leq n$.
	\item $\sum^{n}_{j=1,j\neq k,l}\frac{1}{p_{j}}+\frac{2}{p_{k}}+\frac{2}{p_{l}}\leq n+\frac{1}{r},$ where $1\leq k,l\leq n$ and $k\neq l.$
\end{enumerate} 
\end{theorem}
In order to investigate the sufficiency of conditions listed as above for boundedness of the operator $T$ and to state the corresponding results we need to first describe the necessary region. 

\subsection*{Necessary region $\mathcal R$ and its endpoints} 

Let $\mathcal R=\mathcal R(n)$ denote the closed and bounded region in $\R^{n+1}$ enclosed by the hyperplanes determined by the necessary conditions described in Theorem~\ref{mainresult:nec}. We will refer to it as the necessary region $\mathcal R.$  In order to understand the necessary region $\mathcal R$, we need to find its vertices. The vertices of $\mathcal R$ will be referred to as the endpoints. We will see that in dimension $n=2$ it is easy to write down all the endpoints however for large dimensions the number of endpoints is large and it becomes little difficult to describe all of them. Further note that in view of the multilinear interpolation theory, see~\cite{BS,Ch}, it is enough to prove boundedness of the operator $T$ at the endpoints of the convex region $\mathcal R$. Therefore, knowing the endpoints is important to prove the sufficient part to Theorem~\ref{mainresult:nec} which forms the major part of the paper. Another property that will play a crucial role is the fact that dual of an endpoint, in the sense as described in the previous section, remains an endpoint of $\mathcal R$. This can be easily verified and we skip the detail. We will make use of this fact to identify the endpoints. We will see that at some of the endpoints the operator $T$ fails to satisfy strong type estimates. In this scenario, some boundary points become important provided there holds a strong type result at these points. We shall have positive results at some boundary points in our analysis. For an easy reference we keep the notation same as in~\cite{jbak,oberlin} to denote the points already discussed in there. The endpoints of the region $\mathcal R$ are described as follows.  
\begin{itemize}
\item Clearly $O=(0,0,\dots,0;0)$ is an endpoint.
\item Point $B=(\frac{n-1}{n+1},\frac{n-1}{n+1},\dots,\frac{n-1}{n+1};0)$ as  intersection of $\frac{n+1}{p}=n-1+\frac{2}{r}$ and $\frac{1}{r}=0$.  The dual point is $G=(\frac{n-1}{n+1},\dots,\frac{n-1}{n+1},1;\frac{2}{n+1})$. Note that it is intersection of  $\sum^{n-1}_{j=1}\frac{1}{p_{j}}+\frac{2}{p_{n}}=n-1+\frac{2}{r}$, $\sum^{n-1}_{j=1,j\neq k}\frac{1}{p_{j}}+\frac{2}{p_{k}}+\frac{2}{p_{n}}=n+\frac{1}{r}$ and $p_{n}=1$. There are $n$ different points of this type.

\item Point $M=(\frac{n+1}{n+3},\frac{n+1}{n+3},\dots,\frac{n+1}{n+3};\frac{2}{n+3})$ as intersection of $\frac{n+1}{p}=n-1+\frac{2}{r}$ and $\frac{n+2}{p}=n+\frac{1}{r}$. 
\item Point $E=(\frac{n-1}{n},\frac{n-1}{n},\dots,\frac{n-1}{n},0;0)$ as the intersection of $\sum^{n}_{j=1,j\neq k}\frac{1}{p_{j}}+\frac{2}{p_{k}}= n-1+\frac{2}{r},~k\neq n$, $\frac{1}{p_{n}}=0$ and $\frac{1}{r}=0$. Note that due to symmetry there are $n$ different points of this type. Point $E$ has two different type of dual points. One is of type $P=(\frac{n-1}{n},\frac{n-1}{n},\dots,\frac{n-1}{n},1;1)$ with $n$ different points. The other type of dual is 
$K=(\frac{n-1}{n},\frac{n-1}{n},\dots,\frac{n-1}{n},1,0;\frac{1}{n})$, which can be seen as intersection of $\sum^{n}_{j=1,j\neq n-1}\frac{1}{p_{j}}+\frac{2}{p_{n-1}}=n-1+\frac{2}{r}$, $\sum^{n}_{j=1,j\neq k,n-1}\frac{1}{p_{j}}+\frac{2}{p_{k}}+\frac{2}{p_{n-1}}=n+\frac{1}{r} ,(k\neq n)$, $\frac{1}{p_{n-1}}=1$ and $ \frac{1}{p_{n}}=0$. Note that there are $n(n-1)$ many different points of this type.

\item Point $A=(\frac{n+1}{n+2},\frac{n+1}{n+2},\dots,\frac{n+1}{n+2};1)$ as  intersection of $\frac{n+2}{p}=n+\frac{1}{r}$ and $\frac{1}{r}=1$. Note that $A$ is an endpoint for $n\geq 3$. When $n=2$, the point $A$ lies on the line segment joining $P$ and $P'$ (see Figure~$2$). The dual of $A$ is given by  $A^*=(\frac{n+1}{n+2},\frac{n+1}{n+2},\dots,\frac{n+1}{n+2},0;\frac{1}{n+2})$. There are $n$ different points of this type.

\item It is easy to see the point $C=(0,0,\dots,0,1,0,\dots,0;1)$ is an endpoint and there are $n$ points of this type. 

\item  For $n\geq 3$, the point $Z=(1,1,\dots,1,1,0;1)$ is the intersection of $\sum^{n}_{j=1,j\neq k,l}\frac{1}{p_{j}}+\frac{2}{p_{k}}+\frac{2}{p_{l}}= n+\frac{1}{r},(k,l\neq n)$, $\frac{1}{r}=1$ and $p_{j}=1$, for $j=1,2,\cdots n-1$. There are $n$ different points of this type. The dual point is given by  $Z^*=(1,1,\dots,1,0,0;0)$ which is intersection of  $\sum^{n}_{j=1,j\neq k}\frac{1}{p_{j}}+\frac{2}{p_{k}}= n-1+\frac{2}{r}, k\neq n-1,n$; $\frac{1}{r}=0$ and $p_{j}=1, j=1,2,\dots, n-2$. There are $\frac{n(n-1)}{2}$ points of this type. Note that there is no analogue of $Z$ and $Z^*$ for $n=2.$
\item When $n=3$, consider the point $N=(\frac{3}{5},\frac{3}{5},\frac{1}{5};0)$. It has two type of dual points given by $N^{*1}=(1,\frac{3}{5},\frac{1}{5};\frac{2}{5})$ and $N^{*3}=(\frac{3}{5},\frac{3}{5},1;\frac{4}{5})$. Note that $N$ is not an endpoint but lies on the line segment joining $E$ and $B$. We shall see that $T$ is of strong type at $N$. 
\end{itemize}
\begin{remark}\label{rem}Note that in the above if we interchange positions of $p_j$ for two different values of $j$ we get another endpoint. This is due to the symmetry of the operator $T$. Points obtained in this fashion will be referred to as similar points. We will state results and demonstrate the proofs only for one point of each type and the corresponding results hold for points that are similar to the ones described. 
\end{remark}
The points $O,B,M$ and $A$ have $p_j=p$ for all $j$. These points  have already been addressed in~\cite{oberlin}. Even though the point $N$ is not an endpoint, it plays an important role in proving strong type estimates on some part of the boundary. The region $\mathcal R$ is the closed convex hull of all the endpoints listed above (along with their similar points) in $\R^{n+1}$.

We list down the endpoints for $n=2$ case separately and it is possible to draw the region $\mathcal R$ in this case, see Figure~$2.$ The endpoints where we have strong type results are marked with boldfaced points. For $n=2$ the endpoints are:  $O=(0,0;0), E'=(\frac{1}{2},0;0), E=(0,\frac{1}{2};0), B=(\frac{1}{3},\frac{1}{3};0), K=(0,1;\frac{1}{2}), K'=(1,0;\frac{1}{2}), M=(\frac{3}{5},\frac{3}{5};\frac{2}{5}), C=(0,1;1), C'=(1,0;1), P=(\frac{1}{2},1;1)$, 
$P'=(1,\frac{1}{2};1),  G'=(1,\frac{1}{3};\frac{2}{3}), G=(\frac{1}{3},1;\frac{2}{3}).$ 

\begin{figure} \label{fig1}
	\includegraphics[width=.6\textwidth]{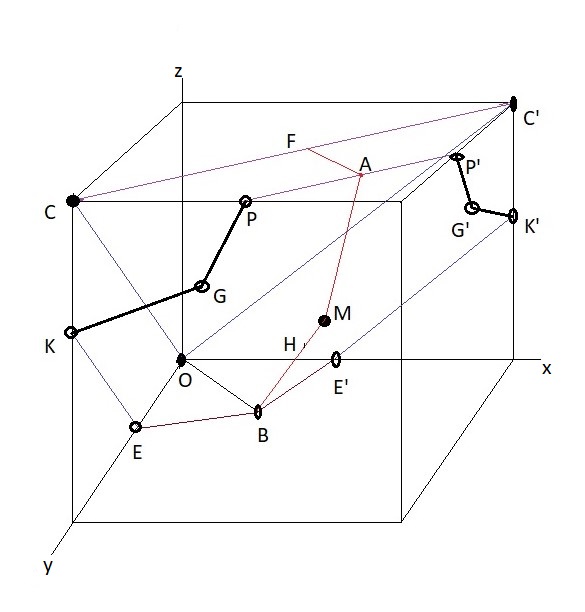}	
	\caption{Region $\mathcal R$ for $n=2$}	
\end{figure}
\subsection{Sufficient part : The case of $n=2$}\label{two}
 Note that for $n=2$ the operator $T$ takes a simpler form as compared to the higher dimensional analogues. Also, there are fewer endpoints in this case and some of the proofs for $n=2$ provide foundation to deal with the case $n\geq 3$. Therefore, we deal with the case of $n=2$ separately. This treatment also helps us understand the problem better. 

As mentioned previously it is enough to prove boundedness of $T$ at the endpoints of $\mathcal R$.  Boundedness of $T$ at endpoints for which all the $p_j's$ are equal is already known due to~\cite{jbak,oberlin}. We include these known points in our statement for completion and provide proofs for the remaining points. 
\begin{theorem}\label{mainresult:sufn=2}{\bf (Sufficient part for $n=2$)}\label{ $n=2$}  In dimension $n=2$ the following estimates hold.
\begin{enumerate}
		\item $T$ is of strong type at $O,C$ and $M$.
		\item $T$ is of restricted type at $P,G,E,$ and $B$. Moreover, $T$ fails to be of strong type at these points.
		\item  $T$ is of weak type at the point $K$ and fails to be of strong type at $K$.
		\item $T$ does not satisfy strong type estimates at points lying on the open line segment $BE.$ 
	\end{enumerate}
\end{theorem}

\begin{theorem}\label{mainresult:sufall}{\bf (Sufficient part for $n\geq3$)} Let $n\geq 3$. Then, 
\begin{enumerate}
	\item  $T$ is of strong type at $O,B,M,A,A^{*},C,Z,Z^*$ and $G$.
	\item $T$ is of restricted type at $E$ and $P$ and it is of restricted weak type at $K$.
	\item Moreover, when $n=3$, $T$ is of strong type at the boundary point $N=(\frac{3}{5},\frac{3}{5},\frac{1}{5};0)$ and its dual points $N^{*1}$ and $N^{*3}$.
\end{enumerate}	
\end{theorem}
\begin{remark} We have the following remarks concerning the results stated as above. 
	\begin{enumerate}
		\item In Theorems~\ref{mainresult:sufn=2} and~\ref{mainresult:sufall} we have stated the result only for one point of each type. The analogous estimates hold at points which are similar to the ones given in theorems. For example, at point $C'$ (and line segment $BE'$) we have the analogous estimate as that of $C$ (and line segment $BE$). 
		\item Since the number of endpoints for $n\geq 3$ is large and different type of points have different type of estimates, we do not record results obtained by applying the standard multilinear interpolation arguments to the estimates obtained in Theorems. We have written down certain positive and negative results separately for points lying on the boundary as they require additional arguments along with multilinear interpolation. For example, strong type estimates on the line segment $PP'$ are not recorded in Theorem~\ref{mainresult:sufn=2} as they follow using the multilinear interpolation Lemma~\ref{christ}. However, the failure of strong type estimates on the line segment $BE$ needs to be discussed through examples and hence we have it in Theorem~\ref{mainresult:sufn=2}.
	\end{enumerate}
\end{remark}
\section{Proof of Theorem~\ref{mainresult:nec}}\label{proof:nec}

	In this section we obtain necessary conditions on exponents for $L^p$ boundedness of the operator $T$. We work with examples considered in \cite{oberlin}. Assume that $T$ is of strong type at $(p_1,p_2,\dots,p_n;r)$ where $1\leq p_1,p_2,\dots,p_n,r\leq \infty$.  
	
	For a positive number $L>0$ consider the following setting. Let $f_{j}=\chi_{_{[-L,L]}},~j=1,2,\dots,n$. Observe that for $|x|<L-1$ we have that $T(f_1,f_2,\dots,f_n)(x)\gtrsim 1$ . Therefore, the assumption on $T$ implies that 
	$$L^{\frac{1}{r}}\lesssim L^{\frac{1}{p_{1}}+\dots+\frac{1}{p_{n}}}$$ 
	for arbitrary large numbers $L$. This yields the first necessary condition in Theorem~\ref{mainresult:nec}, namely  $\frac{1}{r}\leq\sum^{n}_{j=1}\frac{1}{p_{j}}.$ 
	
	Next, consider the functions $f_{j}=\chi_{[-\epsilon,\epsilon]}, j=1,2,\dots,n-1$ and $f_{n}=\chi_{[1-2\epsilon^{2},1+2\epsilon^{2}]}$, where $\epsilon>0$. With this choice of functions for $|x|\leq c\epsilon^{2}$ we have that $T(f_1,f_2,\dots,f_n)(x)\gtrsim \epsilon^{n-1},$ where $c$ is a constant. As earlier we get that $$\epsilon^{n-1+\frac{2}{r}}\lesssim \epsilon^{\frac{1}{p_{1}}+\dots+\frac{1}{p_{n-1}}+\frac{2}{p_{n}}}.$$ Letting $\epsilon\rightarrow 0$, we get the second necessary condition $\frac{1}{p_{1}}+\dots+\frac{1}{p_{n-1}}+\frac{2}{p_{n}}\leq n-1+\frac{2}{r}$ for $k=n.$ Interchanging the roles of functions suitably we get the condition for other values of $k.$
	
	Finally, let $f_{k}=\chi_{[\frac{1}{\sqrt{2}}-\epsilon^{2}, \frac{1}{\sqrt{2}}+\epsilon^{2}]}$, $f_{l}=\chi_{[\frac{-1}{\sqrt{2}}-\epsilon^{2}, \frac{-1}{\sqrt{2}}+\epsilon^{2}]}$ and $f_{j}=\chi_{[-\epsilon,\epsilon]}$ for $j\neq k,l$. Let $J$ denote the box in $\mathbb{R}^{n}$ given by $\chi_J(x_{1},\dots,x_{n})=\prod^{n}_{j=1}f_{j}(x_{j})$. Now observe that the measure of the surface $(J+(t,t,\dots,t))\cap \mathbb{S}^{n-1}$ is of the order of $\epsilon^{n}$ for $|t|\leq \frac{\epsilon}{2n}$. Boundedness of $T$ at the point $(p_{1},p_{2},\dots,p_{n};r)$ implies that $$\epsilon^{n+\frac{1}{r}}\lesssim \epsilon^{\sum^{n}_{j=1,j\neq l,k}\frac{1}{p_{j}}+\frac{2}{p_{k}}+\frac{2}{p_{l}}}.$$ 
	Since $\epsilon$ is arbitrarily small we get the third necessary condition that  $\sum^{n}_{j=1,j\neq k,l}\frac{1}{p_{j}}+\frac{2}{p_{k}}+\frac{2}{p_{l}}\leq n+\frac{1}{r}.$
\qed
\section{Proof of Theorem~\ref{mainresult:sufn=2}}\label{proof:n=2}
In this section we prove Theorem~\ref{mainresult:sufn=2}. We demostrate the arguments at each point listed in the theorem separately. 

Note that in the case of $n=2$ the operator $T$ is given by
\begin{eqnarray*}\label{def}T(f_{1},f_{2})(x)&=&\int_{\mathbb{S}^{1}}f_{1}(x-\sigma_1)f_{2}(x-\sigma_2)d\sigma\\
&=&\int_{0}^{2\pi}f_{1}(x-\cos t)f_{2}(x-\sin t)dt,~x\in \R.
\end{eqnarray*}
Without loss of generality we may assume that $f_1$ and $f_2$ are non-negative functions.

\subsection*{Endpoints $O=(0,0;0), B=(\frac{1}{3},\frac{1}{3};0)$ and  $M=(\frac{3}{5},\frac{3}{5};\frac{2}{5})$} These points are already considered in~\cite{jbak,oberlin}. The operator $T$ satisfies strong type estimates at  $O=(0,0;0)$ and $M=(\frac{3}{5},\frac{3}{5};\frac{2}{5})$. However, it fails to be of strong type at $B=(\frac{1}{3},\frac{1}{3};0).$ We have weaker estimates at this point, namely $T$ is bounded from $L^{3,s}(\R)\times L^{3,s}(\R)\rightarrow L^{\infty}(\R)$ if $0<s\leq 2$. See, Theorem~\ref{jbak1} and \cite{jbak,oberlin} for detail. 

\subsection*{Endpoints $C=(0,1;1)$ and  $P=(\frac{1}{2},1;1)$}
First, note that it is easy to verify that $T$ is of strong type at $C$. Whereas from~ \cite{jbak} the estimate at $P$ is known, namely, $T$ maps $L^{2,1}(\R)\times L^{1}(\R)\rightarrow L^{1}(\R).$ Indeed, one can show that $T$ satisfies strong type estimates at point $(p_1,p_2;r)$ lying on the line segments $CP$, except the point $P$.  

For, note that an arbitrary point on the open line segment $CP$ can be written as $(\frac{1}{q},1;1)$ for $2<q<\infty$. Let $h\in L^{\infty}(\mathbb{R})$\\ and make a change of variables to get that 
\begin{eqnarray*}
	|\langle T(f_1,f_2),h\rangle|&=&|\int_{\mathbb{R}}\left(\int^{1}_{0}f_1(x-\sqrt{1-t^{2}})f_2(x-t)\frac{dt}{\sqrt{1-t^{2}}}\right)h(x)dx|\\
	&=&|\int_{\mathbb{R}}f_2(x)\left(\int^{1}_{0}h(x+t)f_1(x+t-\sqrt{1-t^{2}})\frac{dt}{\sqrt{1-t^{2}}}\right)dx|\\
	&\leq&\Vert f_2\Vert_{L^{1}}\Vert h\Vert_{L^{\infty}(\mathbb{R})}\Vert f_1\Vert_{L^{q}}\Vert \frac{\chi_{[0,1]}(\cdot)}{\sqrt{(1-\cdot)}}\Vert_{L^{q'}}.
\end{eqnarray*}
Since $1<q'<2$, we have $\Vert \frac{\chi_{[0,1]}(\cdot)}{\sqrt{(1-\cdot)}}\Vert_{L^{q'}}<\infty$. This completes the proof. 

Next, we show that $T$  cannot be of strong type at $P$. 
For $h\in L^{\infty}(\mathbb{R})$, we can write
\begin{eqnarray*}
	\langle T(f_1,f_2),h\rangle 
	&=&\int_{\mathbb{R}}f_2(x+1)\left(\int^{1}_{-1}h(x+1+t)f_1(x+1+t-\sqrt{1-t^{2}}) \frac{dt}{\sqrt{1-t^{2}}}\right)dx.
\end{eqnarray*}
Therefore, it is enough to show that  
$$T_{1}(f_1,h)(x)=\int^{1}_{-1}h(x+1+t)f_{1}(x+1+t-\sqrt{1-t^{2}}) \frac{dt}{\sqrt{1-t^{2}}}$$
is a unbounded function for some choice of $h\in L^{\infty}(\mathbb{R})$ and $f_1\in L^{2}(\mathbb{R})$. 

We choose $h=1$ and $f_1(t)=\chi_{[0,\frac{9}{10}]}(|t|)|t|^{-\frac{1}{2}}|\log|t||^{-\frac{2}{3}}.$ Note that $f_1\in L^{2}(\mathbb{R})$. For arbitrarily small $x>0$ we get that,
\begin{eqnarray*}
	T_{1}(f_1,h)(x)&\gtrsim&\int^{1-\frac{x^{2}}{2}}_{\frac{3}{4}}|x+1-t-\sqrt{1-t^{2}}|^{-\frac{1}{2}}|\log|x+1-t-\sqrt{1-t^{2}}||^{-\frac{2}{3}}\frac{dt}{\sqrt{1-t^{2}}}\\
	&\gtrsim&|\log x|^{-\frac{2}{3}}\int^{1-\frac{x^{2}}{2}}_{\frac{3}{4}}(1-t)^{-1}dt.
\end{eqnarray*}
The second inequality in the above estimate follows by using that 
$|x+1-t-\sqrt{1-t^{2}}|\lesssim 1-t$ and $|x+1-t-\sqrt{1-t^{2}}|\gtrsim x^{2}$ for  $t\in [\frac{3}{4},1-\frac{x^{2}}{2}]$.
This gives the desired result.


\subsection*{Endpoint $E=(0,\frac{1}{2};0)$ } We shall show that the operator $T$ does not satisfy of strong type at $E=(0,\frac{1}{2};0).$ However, it maps $L^{\infty}(\R)\times L^{2,1}(\R) $ into $L^{\infty}(\R)$. 

Note that after taking out the $\|f_1\|_{\infty}$, it is enough to show that the operator 
\begin{eqnarray}\label{Ave}
	\tilde{\mathcal{A}}f_{2}(x)=\int_{\mathbb {S}^1}f_{2}(x-\sigma_2)d\sigma,~x\in \R
	\end{eqnarray}
maps $L^{2,1}(\R)$ into $L^{\infty}(\R)$. This follows using H\"{o}lder's inequality. This point is also used in~\cite{jbak}. Infact, the same argument yields strong type estimate for $T$ on the line segment $OE$ except at point $E$.  

Next, we show that the operator $T$ does not verify strong type estimate at $E.$

For, let $f_1=1$ and $f_2(t)=\chi_{[0,\frac{1}{2}]}(|1+t|)|1+t|^{-\frac{1}{2}}|\log|1+t||^{-\frac{2}{3}}$ and note $f_2\in L^2(\R)$.  Let $x\in \R$ be a small negative number, then we have the following.  
\begin{eqnarray*}
T(f_1,f_2)(x)& \gtrsim & \int^{1}_{0}f_2(x-t)\frac{dt}{\sqrt{1-t^{2}}}\\
&\gtrsim & \int^{1-2|x|}_{\frac{1}{2}}\chi_{[0,\frac{1}{2}]}(|1-t+x|)|1-t+x|^{-\frac{1}{2}}|\log|1-t+x||^{-\frac{2}{3}}\frac{dt}{\sqrt{1-t^{2}}}.
\end{eqnarray*}
Note that $|1+x-t|\lesssim |1-t|$ and $|1+x-t|\gtrsim |x|$. Therefore, we get that 
\begin{eqnarray*}
T(f_{1},f_{2})(x)&\gtrsim& \int^{1-2|x|}_{\frac{1}{2}}|\log|x||^{-\frac{2}{3}}|1-t|^{-1}dt\\
&=&|\log|x||^{-\frac{2}{3}}\int^{1-2|x|}_{\frac{1}{2}}|1-t|^{-1}dt\\
&\gtrsim &|\log|x||^{-\frac{2}{3}}|\log2|x||.
\end{eqnarray*}
Clearly, the function in the estimate above is not bounded near the origin $x=0$. 

\subsection*{Endpoint $K=(0,1;\frac{1}{2})$} The operator $T$ is of weak type at $K$, i.e. it maps $L^{\infty}(\mathbb{R})\times L^{1}(\mathbb{R}) \rightarrow L^{2,\infty}(\mathbb{R})$. Moreover, $T$ fails to be of strong type at $K$. 

Note that in view of the standard duality arguments (see \cite{Grafakosclassical}, page $69$), the weak type estimates for the operator $T$ can be deduced by considering the following estimate.
\begin{eqnarray*}
\sup_{\Vert h\Vert_{L^{2,1}}\leq 1}|\int_{\mathbb{R}}T(f_1,f_2)(x)h(x)dx|\lesssim \Vert f_1\Vert_{L^{\infty}} \Vert f_2\Vert_{L^{1}}.
\end{eqnarray*}
Consider
	\begin{eqnarray*}
	|\int_{\mathbb{R}}T(f_1,f_2)(x)h(x)dx|&=&|\int^{\infty}_{-\infty}\int^{1}_{-1}f_2(x-t)f_1(x-\sqrt{1-t^{2}})\frac{dt}{\sqrt{1-t^{2}}}h(x)dx|\\
	&=&|\int^{\infty}_{-\infty}f_2(x)\int^{1}_{-1}h(x+t)f_1(x+t-\sqrt{1-t^{2}})\frac{dt}{\sqrt{1-t^{2}}}dx|\\
	&\lesssim&\Vert f_2\Vert_{L^{1}}\Vert f_1\Vert_{L^{\infty}}\Vert h\Vert_{L^{2,1}}.
\end{eqnarray*}
Next, we see that the operator $T$ cannot be bounded from $ L^{\infty}(\mathbb{R})\times L^{1}(\mathbb{R})\rightarrow L^{2}(\mathbb{R})$. 

Set $f_1=1$ and consider 
\begin{eqnarray*}
T(f_1,f_2)(x)&=&\int_{\mathbb{S}^{1}}f_2(x-z)d\sigma(y,z)\\
&=&\int_{\mathbb{R}}\hat{f}_2(\xi)e^{2\pi\iota x\cdot\xi}\left(\int_{\mathbb{S}^{1}}e^{-2\pi\iota(0,\xi)\cdot(y,z)}d\sigma(y,z)\right)d\xi\\
&=&\int_{\mathbb{R}}\hat{f}_2(\xi)\widehat{d\sigma}(0,\xi)e^{2\pi\iota x\cdot\xi}d\xi.
\end{eqnarray*}
Therefore, we get that 
\begin{eqnarray}\label{eq101}
\Vert T(f_1,f_2)\Vert^{2}_{L^{2}}=\Vert \widehat{T(f_1,f_2)}\Vert^{2}_{L^{2}}\gtrsim \int_{|\xi|>1}|\hat{f_{2}}(\xi)|^{2}(1+|\xi|)^{-1}d\xi.
\end{eqnarray}
Here we have used decay estimate for the Fourier transform of surface measure, namely,  
$$\widehat{d\sigma}(0,\xi)\sim C |\xi|^{\frac{-1}{2}}[e^{2\pi\iota|\xi|}\sum^{\infty}_{j=0}\alpha_{j}|\xi|^{-j}+e^{-2\pi\iota|\xi|}\sum^{\infty}_{j=0}\beta_{j}|\xi|^{-j}],$$ 
as $|\xi|\rightarrow\infty$, for suitable constants $\alpha_{j},\beta_{j}$ (see \cite{Steinbook}, page $391$). 
Choose $f_2\in L^{1}(\R)$ such that $|\hat f_{2}(\xi)|$ decays slower than $(\log|\xi|)^{-\frac{1}{2}}$ for $|\xi|>1.$ This implies that integral in the estimate ~\ref{eq101} diverges and consequently we get the desired result. 

\subsection*{Endpoint $G=(\frac{1}{3},1;\frac{2}{3})$}
 We show that the operator $T$ maps $ L^{3,\frac{3}{2}}(\mathbb{R})\times L^{1}(\mathbb{R})\rightarrow L^{\frac{3}{2}}(\mathbb{R}).$ 
  Subsequently, we get that $T$ is of restricted type at $G$. 
Observe that for the said boundedness result, it is enough to prove that 
\begin{eqnarray*}
|\langle T(f_1,f_2),h\rangle|\lesssim \Vert f_1\Vert_{L^{3,\frac{3}{2}}}\Vert f_2\Vert_{L^{1}}\Vert h\Vert_{L^{3}},~~~ h\in L^{3}(\mathbb{R}).
\end{eqnarray*}
Consider
\begin{eqnarray*}
\langle T(f_1,f_2),h\rangle&=&\int_{\mathbb{R}}\left(\int^{2\pi}_{0}f_1(x-\cos\theta)f_2(x-\sin\theta)d\theta\right)h(x)dx\\
&=&\int_{\mathbb{R}}f_2(x)\left(\int^{2\pi}_{0}f_1(x+\sin\theta-\cos\theta)h(x+\sin\theta)d\theta\right)dx\\
&=&\int_{\mathbb{R}}f_2(x)\left(\int_{\mathbb{S}^{1}}f_1(x-\sigma\cdot v_{1})h(x-\sigma\cdot v_{2})d\sigma\right)dx,	
\end{eqnarray*}
where $v_{1}=\left( \begin{array}{c} -1 \\ 1 \end{array} \right)$ and $v_{2}=\left( \begin{array}{c} 0 \\ 1 \end{array} \right)$ are linearly independent vectors in $\mathbb{R}^{2}$ and $\sigma=e^{i\theta}\in\mathbb{S}^{1}$. 
We consider the following operator
\begin{eqnarray*}
T_v(f_{1},h)(x)=\int_{\mathbb{S}^{1}}f_1(x-\sigma\cdot v_{1})h(x-\sigma\cdot v_{2})d\sigma.
\end{eqnarray*}
Using the same argument as in case of point $E$, we have the following 
\begin{eqnarray}\label{T_v}
\Vert T_v(f_{1},h)\Vert_{L^{\infty}}\lesssim \Vert h\Vert_{L^{2,1}} \Vert f_{1}\Vert_{L^{\infty}}~~\text{and}~~\Vert T_v(f_{1},h)\Vert_{L^{\infty}}\lesssim  \Vert h\Vert_{L^{\infty}}\Vert f_{1}\Vert_{L^{2,1}}.
\end{eqnarray} 
For a small positive number $0<\epsilon<\frac{1}{4}$, we decompose the operator as $T_v(f_{1},h)(x)\lesssim \sum_{j=1,2}U_{j}(f_{1},h)(x)$ for almost every $x\in \mathbb{R}$, where  
\begin{eqnarray*}
U_{j}(f_{1},h)(x)=\int_{\{\sigma\in\mathbb{S}^{1}:|\sigma\cdot v_{j}|>\epsilon\}}f_1(x-\sigma\cdot v_{1})h(x-\sigma\cdot v_{2})d\sigma, ~~j=1,2.
\end{eqnarray*}
It is easy to show that (also see Lemma 2 in~\cite{jbak})
\begin{eqnarray}\label{U_j}
\Vert U_{1}(f_{1},h)\Vert_{L^{\infty}}\lesssim \Vert h\Vert_{L^{\infty}}\Vert f_{1}\Vert_{L^{1}}~~\text{and}~~
\Vert U_{2}(h,f_{1})\Vert_{L^{\infty}}\lesssim \Vert h\Vert_{L^{1}}\Vert f_{2}\Vert_{L^{\infty}}
\end{eqnarray} 
Observe that $U_j$ also satisfies the estimate~(\ref{T_v}). Therefore, interpolating between these estimates for $U_j$ we get that the desired result holds for $U_j$. Consequently, we get that 
$$\Vert T_v(f_{1},h)\Vert_{L^{\infty}}\lesssim \Vert h\Vert_{L^{3}}\Vert f_{1}\Vert_{L^{3,\frac{3}{2}}}.$$ 
This yields that 
\begin{eqnarray*}
|\langle T(f_1,f_2),h\rangle| &\leq& \Vert f_2\Vert_{L^{1}}\Vert T_v(f_1,h)\Vert_{L^{\infty}}\\
&\lesssim &\Vert f_2\Vert_{L^{1}}\Vert f_{1}\Vert_{L^{3,\frac{3}{2}}}\Vert h\Vert_{L^{3}}.
\end{eqnarray*}

Next,  we give an example to show that $T$ fails to be of strong type at $G$, i.e., it is unbounded from $L^{3}(\mathbb{R})\times  L^{1}(\mathbb{R})$ into $L^{\frac{3}{2}}(\mathbb{R})$.

For, let $f_2\in L^1(\R)$ and $f_1, h\in L^{3}(\mathbb{R})$ and consider  
\begin{eqnarray*}
\langle T(f_1,f_2),h\rangle &=& \int_{\mathbb{R}}f_1(x) T_v(f_1,h)(x)dx	\end{eqnarray*}
where $T_{{v}}(f_1,h)(x)=\int^{1}_{0}h(x+t)f_{1}(x+t-\sqrt{1-t^{2}})\frac{dt}{\sqrt{1-t^{2}}}$ is the same operator as previously.  Since boundedness properties of the operators $T_v$ and $T$ are equivalent, it suffices to prove that $T_{{v}}(f_1,h) \notin L^{\infty}(\mathbb{R})$ for a suitable choice of functions $h$ and $f_1$ in $L^{3}(\mathbb{R})$. 

Let $f_1(t)=\chi_{[0,\frac{9}{10}]}(|t|)|t|^{-\frac{1}{3}}|\log|t||^{-\frac{2}{5}}$ and $h(t)=f_{1}(1+t)$. For a small negative real number $x$, one has the following. 
\begin{eqnarray*}
T_v(f_1,h)(x)\geq\int^{1-2|x|}_{\frac{1}{2}}\frac{|1+x-t|^{-\frac{1}{3}}}{|\log|1+x+t||^{\frac{2}{5}}}\frac{|x-\sqrt{1-t^{2}}|^{-\frac{1}{3}}}{|\log|x-\sqrt{1-t^{2}}||^{\frac{2}{5}}}\frac{dt}{\sqrt{1-t^{2}}}
	\end{eqnarray*}
Observe that for our choice of $x$ and $t$ in the integral above we have $|1+x-t|\lesssim |1-t|$, $|x-\sqrt{1-t^{2}}|\lesssim \sqrt{(1-t)}$, $|\log|1+x-t||\lesssim |\log|x||$ and $|\log|x-\sqrt{1-t^{2}}||\lesssim |\log|x||$. Therefore, 
\begin{eqnarray*}
T_v(f_1,h)(x)&\gtrsim& |\log|x||^{-\frac{4}{5}}\int^{1-2|x|}_{\frac{1}{2}}|1-t|^{-1}dt\\
&=&|\log|x||^{-\frac{4}{5}}(-\log2|x|-\log2).
\end{eqnarray*}
This yields the desired result.

\subsection*{Open line segment $BE$} Let $(\frac{1}{p_1},\frac{1}{p_2};0)$ be a point on the open line segment $BE$. Note that one can write $\frac{1}{p_{1}}=\frac{\theta}{3}$ and $\frac{1}{p_{2}}=\frac{\theta}{3}+\frac{1-\theta}{2}$ for $\theta\in (0,1)$. Set $\epsilon=\frac{1-\theta}{6}$ and write   $\frac{1}{p_{1}}=\frac{1}{3}-2\epsilon$ and $\frac{1}{p_{2}}=\frac{1}{3}+\epsilon$. 
Consider the functions  $f_{1}(t)=\chi_{[0,\frac{9}{10}]}(|t|)|t|^{-(\frac{1}{3}-2\epsilon)}|\log|t||^{-\frac{1}{3}}$ and  $f_{2}(t)=\chi_{[0,\frac{9}{10}]}(|1+t|)|1+t|^{-(\frac{1}{3}+\epsilon)}|\log|1+t||^{-\frac{1}{2}}$. Then for $x$ near the origin we have 
\begin{eqnarray*}
	T(f_{1},f_{2})(x)&\geq&\int^{1}_{0}f_{1}(x-\sqrt{1-t^{2}})f_{2}(x-t)\frac{dt}{\sqrt{1-t^{2}}}\\
	&\gtrsim& \int^{1-2|x|}_{\frac{3}{4}}f_{1}(x-\sqrt{1-t^{2}})f_{2}(x-t)\frac{dt}{\sqrt{1-t^{2}}}.
	\end{eqnarray*}
	Observe that for arbitrarily small $x$, $|1+x-t|\lesssim |1-t|$, $|1+x-t|\gtrsim |x|$, $|x-\sqrt{1-t^{2}}|\lesssim \sqrt{1-t^{2}}$ and $|x-\sqrt{1-t^{2}}|\gtrsim |x|$.
	Therefore we get,
	\begin{eqnarray*}
	T(f_{1},f_{2})(x)&\gtrsim& |\log|x||^{-\frac{5}{6}}\int^{1-2|x|}_{\frac{3}{4}}|1-t|^{-1}dt\\
	&\gtrsim& -|\log|x||^{-\frac{5}{6}}\log|x|.
	\end{eqnarray*}
The above tends to infinity as $x\rightarrow 0$.

This completes the proof of Theorem~\ref{mainresult:sufn=2}.\qed
\begin{remark}
 \begin{enumerate}
 	\item The operator $T$ satisfies strong type estimates at points lying in regions $OEKC$, except at points $E$ and $K$, see Figure~$2.$ Observe that strong type estimates on $OC$ follow by the Riesz-Thorin interpolation. 
 	Next, note that on $OEKC$, we have  $\frac{1}{p_{1}}=0$. Therefore, we have 
 	\begin{eqnarray*}
 		|T(f_{1},f_{2})(x)|\leq \Vert f_{1}\Vert_{L^{\infty}}\tilde{\mathcal{A}}f_{2}(x),
 	\end{eqnarray*}
 where $\tilde{\mathcal A}$ is same as defined earlier in~\ref{Ave} and can be written as 
 	\begin{eqnarray*}
 		\tilde{\mathcal{A}}f_{2}(x)=\int_{\mathbb{R}}\hat{f_{2}}(\xi)\hat{d\sigma}(\xi,0)e^{2\pi \iota x\cdot\xi}d\xi.
 	\end{eqnarray*}
 	Using the estimate $|\hat{d\sigma}(\xi,0)|\lesssim (1+|\xi|)^{-\frac{1}{2}},$ we get the following (see  \cite{Hormander})
 	\begin{eqnarray*}
 	\tilde{\mathcal{A}}: L^{p_{2}}(\mathbb{R})\rightarrow L^{r}(\mathbb{R}),~~\text{for}~~1<p_{2}\leq2\leq r<\infty~~with~~\frac{1}{p_{2}}-\frac{1}{r}\leq\frac{1}{2}. 
 	\end{eqnarray*} 
 This implies strong type estimates for $T$ in the  region $OEKC$, except on the line segments $OE$ and $KC$. The required estimates are already proved for points $OE,$ whereas for points on $KC$ they can be deduced using  the Marcinkiewicz interpolation theorem. 
 \end{enumerate}
\end{remark}
\section{Proof of Theorem~\ref{mainresult:sufall}}\label{proof:all}
 We deal with each point separately. We repeat that we will describe proofs for one endpoint of each type. We use the idea from~\cite{jbak,oberlin}.

\subsection*{Endpoints $O, Z, Z^*, B, M, A, A^*, $ and $ C$}
The boundedness of $T$ at the points $O, Z^*, B, M$ and $A$ is already known due to \cite{jbak,oberlin}. 

Moreover, using standard duality arguments one can deduce corresponding estimates at their dual points. In particular,  boundedness of $T$ at $Z=(1,1,\dots,1,0;1)$ and $A^*=(\frac{n+1}{n+2},\frac{n+1}{n+2},\dots,\frac{n+1}{n+2},0;\frac{1}{n+2})$ can be deduced from that of $Z^*=(1,1,\dots,1,0,0;0)$ and $A=(\frac{n+1}{n+2},\frac{n+1}{n+2},\dots,\frac{n+1}{n+2};1)$ respectively. 

Further, the desired estimates for $T$ at $C$ follows in a straightforward manner.  
\subsection*{Endpoint $E=(\frac{n-1}{n},\dots,\frac{n-1}{n},0;0)$}
We show that the operator $T$ is of restricted type at the point $E$. 

Let $v_{1},v_{2},\cdots ,v_{n}$ be linearly independent vectors in  $\R^n$and consider the operator $T_v$ defined as earlier. Further, let $\Lambda:\mathbb{R}^{n}\rightarrow \mathbb{R}^{n}$ be a linear map from $\R^n$ to $\R^n$ such that $e_{j}\cdot \Lambda x=v_{j}\cdot x,~x\in\mathbb{R}^{n}$. Fix a unit vector $u_{n}$ with $\Lambda u_{n}=c(0,0,\dots,0,1)$ for some  $c\in\mathbb{R}\setminus\{0\}$ and let $\{u_{1},u_{2},\dots,u_{n}\}$ be an orthonormal basis of $\mathbb{R}^{n}$. 

Let $\eta=(\eta_{1},\dots,\eta_{n-1})$ denote an element of $ \mathbb{S}^{n-2}$ and $d\eta$ be the normalised Lebesgue measure on $\mathbb{S}^{n-2}$. We consider the parametrization of $\mathbb{S}^{n-1}$ given by $$\sigma=\sum^{n-1}_{j=1}r\eta_{j}u_{j}+sgn(r)\sqrt{1-r^{2}}u_{n}.$$ where $-1\leq r\leq 1$.

For convenience we will use the notation $\vec{f}(y)=\prod_{j=1}^nf_j(y_j),$ where $y=(y_1,y_2,\dots,y_n)\in \R^n.$ With this we have the following.
\begin{eqnarray*}
	&&T(f_{1},f_{2},\dots,f_{n})(x)\\
	&&=\int_{\mathbb{S}^{n-1}}\vec{f}\big((x,x,\dots,x)-\Lambda\sigma\big)d\sigma\\
	&&=\int^{1}_{-1}\int_{\mathbb{S}^{n-2}}\vec{f}\big((x,x,\dots,x)-sgn(r)\sqrt{1-r^{2}}(0,0,\dots,0,c)-\Lambda(r\sum^{n-1}_{k=1}\eta_{k}u_{k})\big)d\eta |r|^{n-2}\frac{dr}{\sqrt{1-r^{2}}}\\
	&&=\int^{1}_{-1}G(r)|r|^{n-2}\frac{dr}{\sqrt{1-r^{2}}},
\end{eqnarray*}
where $G(r)=\int_{\mathbb{S}^{n-2}}\prod^{n-1}_{j=1}f_{j}\big(x- re_{j}\cdot\Lambda(\sum^{n-1}_{k=1}\eta_{k}u_{k})\big)f_{n}\big(x- sgn(r)\sqrt{1-r^{2}}c-re_{n}\cdot\Lambda(\sum^{n-1}_{k=1}\eta_{k}u_{k})\big)d\eta$. Consider
\begin{eqnarray*}
	&&\int^{1}_{-1}G(r)|r|^{n-2}dr\\
	&&=\int^{1}_{-1} \int_{\mathbb{S}^{n-2}}\prod^{n-1}_{j=1}f_{j}\big(x- re_{j}\cdot\Lambda(\sum^{n-1}_{k=1}\eta_{k}u_{k})\big)f_{n}\big(x- sgn(r)\sqrt{1-r^{2}}c-re_{n}\cdot\Lambda(\sum^{n-1}_{k=1}\eta_{k}u_{k})\big)d\eta|r|^{n-2}dr\\
	&&\lesssim \Vert f_{n}\Vert_{L^{\infty}}\int^{1}_{0}\int_{\mathbb{S}^{n-2}}\prod^{n-1}_{j=1}f_{j}\big(x- re_{j}\cdot\Lambda(\sum^{n-1}_{k=1}\eta_{k}u_{k})\big)d\eta r^{n-2}dr\\
	&&\lesssim \Vert f_{n}\Vert_{L^{\infty}}  \int_{\mathbb{R}^{n-1}}\prod^{n-1}_{j=1}f_{j}\big(x-e_{j}\cdot \Lambda(\sum^{n-1}_{k=1}x_{k}u_{k})\big)dx_{1}dx_{2}\dots dx_{n-1}\\
	&&\lesssim |\Lambda|^{-1}\Vert f_{n}\Vert_{L^{\infty}}\prod^{n-1}_{j=1}\Vert f_{j}\Vert_{L^{1}}.
\end{eqnarray*}

On the other hand we have that 
\begin{eqnarray}
	G(r)&=&\nonumber \int_{\mathbb{S}^{n-2}}\prod^{n-1}_{j=1}f_{j}(x-re_{j}\cdot\Lambda(\sum^{n-1}_{k=1}\eta_{k}u_{k}))f_{n}(x-sgn(r)\sqrt{1-r^{2}}c-re_{n}\cdot \Lambda(\sum^{n-1}_{k=1}\eta_{k}u_{k}))d\eta\\
	&\leq&\nonumber \Vert f_{n}\Vert_{L^{\infty}}\int_{\mathbb{S}^{n-2}}\prod^{n-1}_{j=1}f^{r}_{j}(\frac{x}{r}-e_{j}\cdot\Lambda(\sum^{n-1}_{k=1}\eta_{k}u_{k}))d\eta \\ 
	&\leq& \nonumber r^{-\frac{(n-1)(n-2)}{n}}\Vert f_{n}\Vert_{L^{\infty}}\prod^{n-1}_{j=1}\Vert f_{j}\Vert_{\frac{n}{n-2},1} \\
	&\simeq &  \label{b} r^{-\frac{(n-1)(n-2)}{n}}\Vert f_{n}\Vert_{L^{\infty}}\prod^{n-1}_{j=1} |I_{j}|^{\frac{n-2}{n}},
\end{eqnarray} 
where $f_{j}=\chi_{I_{j}},~j=1,2,\dots,n-1$.

Note that in the above we have used the boundedness of the operator $T$ at the point $B$ in dimension $n-1$. 

Next, we need to consider two cases separately to complete the proof in the following fashion. 

\noindent
\textbf{Case 1:} When $\prod^{n-1}_{j=1}|I_{j}|\geq 1$. In this case we use the estimate~\ref{b} as follows. 
\begin{eqnarray*}
\int^{1}_{-1}G(r)|r|^{n-2}\frac{dr}{\sqrt{1-r^{2}}} 
&\leq &\Vert f_{n}\Vert_{L^{\infty}}\prod^{n-1}_{j=1} |I_{j}|^{\frac{n-2}{n}} \int^{1}_{-1}r^{-\frac{(n-1)(n-2)}{n}}|r|^{n-2}\frac{dr}{\sqrt{1-r^{2}}}\\
&\lesssim& \Vert f_{n}\Vert_{L^{\infty}}\prod^{n-1}_{j=1} |I_{j}|^{\frac{n-2}{n}}\\
&\leq & \Vert f_{n}\Vert_{L^{\infty}}\prod^{n-1}_{j=1} |I_{j}|^{\frac{n-1}{n}}.
\end{eqnarray*} 
\noindent
\textbf{Case 2:} When $\prod^{n-1}_{j=1}|I_{j}|<1$. We choose $\delta=(\prod^{n-1}_{j=1}|I_{j}|)^{\frac{2}{n}}$ and consider
\begin{eqnarray*} 
\int^{1-\delta}_{-1+\delta}G(r)|r|^{n-2}\frac{dr}{\sqrt{1-r^{2}}}
&\leq & \delta^{-\frac{1}{2}}C\int^{1-\delta}_{-1+\delta}G(r)|r|^{n-2}dr\\
	&\lesssim & \delta^{-\frac{1}{2}}\Vert f_{n}\Vert_{L^{\infty}}\prod^{n-1}_{j=1}|I_{j}|\\
	&=&\Vert f_{n}\Vert_{L^{\infty}}\prod^{n-1}_{j=1}|I_{j}|^{\frac{n-1}{n}}.
\end{eqnarray*}
For the other part, we have 
\begin{eqnarray*}
\int_{\{1-\delta\leq |r|\leq 1\}}G(r)|r|^{n-2}\frac{dr}{\sqrt{1-r^{2}}}
&\lesssim& \Vert f_{n}\Vert_{L^{\infty}}\prod^{n-1}_{j=1}|I_{j}|^{\frac{n-2}{n}}\int_{\{1-\delta\leq |r|
\leq  1\}}|r|^{(n-2)-\frac{(n-1)(n-2)}{n}}\frac{dr}{\sqrt{1-r^{2}}}\\
	&\lesssim &  \Vert f_{n}\Vert_{L^{\infty}}\prod^{n-1}_{j=1}|I_{j}|^{\frac{n-2}{n}} \delta^{\frac{1}{2}}\\
	&=&\Vert f_{n}\Vert_{L^{\infty}}\prod^{n-1}_{j=1}|I_{j}|^{\frac{n-1}{n}}.
\end{eqnarray*}
Here in the estimate above we have used (\ref{b}). 


\subsection*{Endpoints $K=(\frac{n-1}{n},\frac{n-1}{n},\frac{n-1}{n},\dots,\frac{n-1}{n},1,0;\frac{1}{n}), P=(\frac{n-1}{n},\frac{n-1}{n},\dots,\frac{n-1}{n},1;1)$ and $G=(\frac{n-1}{n+1},\dots,\frac{n-1}{n+1},1;\frac{2}{n+1})$} First we show that the operator $T_{v}$ (or $T$) is of restricted  weak type at $K$ and of restricted type at  $P$. This can be proved using the boundedness of $T_{v}$ at the point $E$ along with the duality argument. 

In order to prove the boundedness at the point $K$, we need to show that $T_{v}$ maps  $L^{\frac{n}{n-1},1}(\mathbb R)\times\dots\times L^{\frac{n}{n-1},1}(\mathbb R)\times L^{1}(\mathbb R)\times L^{\infty}(\mathbb R)\rightarrow L^{n,\infty} (\mathbb R)$.
It suffices to show that  $$\sup_{\Vert h\Vert_{L^{\frac{n}{n-1},1}}=1}|\langle T_{v}(f_{1},f_{2},\dots,f_{n}),h\rangle|\lesssim \Vert f_{n-1}\Vert_{L^{1}}\Vert f_{n}\Vert_{L^{\infty}}\prod^{n-2}_{j=1}\Vert f_{j}\Vert_{L^{\frac{n}{n-1},1}}.$$
Now, 
\begin{eqnarray*}
|\langle T_{v}(f_{1},f_{2},\dots,f_{n}),h\rangle|
&=&|\langle f_{n-1},T^{*n-1}_{v}(f_{1},f_{2},\dots,f_{n-2},h,f_{n})|\\
&\leq& \Vert f_{n-1}\Vert_{L^{1}}\Vert T^{*n-1}_{v}(f_{1},f_{2},\dots,f_{n-2},h,f_{n})\Vert_{L^{\infty}}. 
\end{eqnarray*}  
Invoking the boundedness of $T_{v}$ at the point $E$ we get the desired estimate. 
\begin{eqnarray*}
|\langle T_{v}(f_{1},f_{2},\dots,f_{n}),h\rangle|\lesssim \Vert f_{n-1}\Vert_{L^{1}} \Vert f_{n}\Vert_{L^{\infty}}\prod^{n-2}_{j=1}\Vert f_{j}\Vert_{L^{\frac{n}{n-1},1}}\Vert h\Vert_{L^{\frac{n}{n-1},1}}.
\end{eqnarray*}
Now, in order to prove restricted type  boundedness at the point $P$, we need to show that  $T_{v}$ maps $L^{\frac{n}{n-1},1}(\mathbb R)\times \dots\times L^{\frac{n}{n-1},1}(\mathbb R)\times L^{1}(\mathbb R)\rightarrow L^{1} (\mathbb R)$. It suffices to show that 
\begin{eqnarray*}
\sup_{\Vert h\Vert_{L^{\infty}}=1}|\langle T_{v}(f_{1},\dots,f_{n}),h\rangle|\lesssim \Vert f_{n}\Vert_{L^{1}}\prod^{n-1}_{j=1}\Vert f_{j}\Vert_{\frac{n}{n-1},1}.
\end{eqnarray*}
This can be proved in a similar manner as the case of boundedness at $K$.
Next, note that the dual of $L^{\frac{n+1}{2}}(\R)$ is $L^{\frac{n+1}{n-1}}(\R)$. Then, using the same reasoning as above, this time with the point $B=(\frac{n-1}{n+1},\dots,\frac{n-1}{n+1};0)$, we get that $T_{v}$ is bounded at $G$. 

\subsection*{Strong type estimates for $T$ at  $N=(\frac{3}{5},\frac{3}{5},\frac{1}{5};0)$} 
When $n=3$ we have the strong typeness at  $N=(\frac{3}{5},\frac{3}{5},\frac{1}{5};0)$. This point lies on the segment joining the points $E=(\frac{2}{3},\frac{2}{3},0;0)$ and $B=(\frac{1}{2},\frac{1}{2},\frac{1}{2};0)$. Note that, we do not have strong type estimates at $E$. Therefore, proving strong type estimates for $T$ at $N$ would give us the same on the boudary between $N$ and $B$. Let $\epsilon>0$ be a small number and consider the following operators (see~\cite{jbak} for more details). 
\begin{eqnarray*}
U_{j}(f_1,f_2,f_3)(x):=\int_{\{\sigma\in\mathbb{S}^{2}:|v_{j}\cdot \sigma|>\epsilon\}}f_1(x-v_{1}\cdot \sigma)f_{2}(x-v_{2}\cdot \sigma)f_3(x-v_{3}\cdot \sigma)d\sigma, ~~j=1,2,3.
\end{eqnarray*} 
See~\cite{jbak} for more detail about $U_{1},U_{2}$ and $U_{3}$. 
We know that $U_{1}$ is bounded at the point $(0,1,1;0)$. Also, it is of restricted weak type  at $E=(\frac{2}{3},\frac{2}{3},0;0)$ and of strong type at $Z^*=(1,0,0;0)$ using the corresponding estimates for the operator $T$. The interpolation result from~\cite{Ch} yields that $U_1$ maps $L^{\frac{5}{3},2}(\mathbb R)\times L^{\frac{5}{3},2}(\mathbb R)\times L^{5,\infty}(\mathbb R)\rightarrow L^{\infty}(\mathbb R)$. Subsequently, we get that $U_{1}$ is of strong type at $N$. In a simlar way, we can prove the strong typeness of $U_2$ at $N$. 

Next, note that $U_{3}$ is of strong type at $(1,1,0;0)$. Further, we know that $T_v$ is strong type bounded at $(0,0,\frac{1}{2};0).$The Riesz-Thorin interpolation theorem for multilinear operators~\cite{BS} yields that $U_3$ is of strong type at $N$. This completes the proof. 

The standard duality arguments imply strong type boundedness of $T$ at dual points of $N$.

This completes the proof of Theorem~\ref{mainresult:sufall}. \qed

\section*{Acknowledgement}The first author acknowledges the financial support from the Science and Engineering Research Board (SERB), Government of India, under the grant MATRICS: MTR/2017/000039/Math. The second author is supported by CSIR (NET), file no. 09/1020 (0094)/2016-EMR-I. 

\end{document}